\documentclass[graybox]{svmult}

 \usepackage{amsmath}
 \usepackage{amssymb}

\usepackage{mathptmx}     
\usepackage{helvet}           
\usepackage{courier}         
\usepackage{type1cm}        
\usepackage{makeidx}         
\usepackage{graphicx}        
\usepackage{multicol}        
\usepackage[bottom]{footmisc}

\input xypic
 \xyoption{all}
\usepackage{hyperref} 

\newcommand{\C}{\mathbb C}
\newcommand{\R}{\mathbb R}
\newcommand{\Q}{\mathbb Q}
\newcommand{\Z}{\mathbb Z}
\newcommand{\F}{\mathbb F}

\def\O{{\cal O}}
\def\H{{\mathbb H}}
\def\a{{\mathfrak a}}
\def\b{{\mathfrak b}}
\def\p{{\mathfrak p}}

\def\gl{{\rm GL}_2}
\def\sl{{\rm SL}_2}
\def\pgl{{\rm PGL}_2}
\def\psl{{\rm PSL}_2}

\makeindex       

\begin{document}

\title*{Arithmetic Aspects of Bianchi Groups}
\author{Mehmet Haluk \c{S}eng\"un}

\institute{ \at Max Planck Institute for Mathematics, Bonn, Germany, \email{haluk@mpi-bonn.mpg.de}
\and Current Address \at University of Warwick, Coventry, UK, \email{M.H.Sengun@warwick.ac.uk}}

\maketitle

\abstract{We discuss several arithmetic aspects of Bianchi groups, especially from a computational point of view. In particular, we consider 
computing the homology of Bianchi groups together with the Hecke action, connections with automorphic forms, abelian varieties, Galois 
representations and the torsion in the homology of Bianchi groups. Along the way, we list several open problems and conjectures, 
survey the related literature, presenting concrete examples and numerical data.}

\section{Introduction}
Let $L$ be an algebraic number field with ring of integers $\O_L$ and signature $(r,s)$. The associated {\em modular group} $\psl(\O_L)$ is a discrete subgroup of the real Lie group $\psl(\R)^r \times \psl(\C)^s$ and thus acts via isometries on the product space $\mathcal{H}_2^r \times \mathcal{H}_3^s$, where $\mathcal{H}_2, \mathcal{H}_3$ denote the hyperbolic spaces of dimensions 2 and 3 respectively.  Any torsion-free finite-index subgroup 
$\Gamma$ of $\psl(\O_L)$ acts freely and properly discontinuously on $\mathcal{H}_2^r \times \mathcal{H}_3^s$ and the quotient manifolds are locally symmetric spaces which are important for the theory of automorphic forms, differential geometry, number theory and representation theory.

When $L$ is totally real (that is, $s=0$), these quotient manifolds can be viewed as the complex points of quasi-projective (complex) algebraic varieties. 
This allows one to employ tools of algebraic geometry in their study and is essential for our
understanding of the web of conjectures within Langlands programme that surrounds $\gl(\O_L)$. 
However in the case $s>0$, due to the fact that $\mathcal{H}_3$ does not admit any complex structure, there is no visible link to (complex) 
algebraic geometry. 

The simplest case beyond totally real fields is the case of imaginary quadratic fields, that is $(r,s)=(0,1)$. The modular groups in this case are called {\bf Bianchi groups} 
after Luigi Bianchi who studied their geometry as early as 1892 \cite{bianchi}. To the best of our knowledge, the earliest modern considerations of Bianchi groups 
go back to the 1960's, notably to \cite{kubota-64, serre-70, swan-68}. The investigations on Bianchi groups related to aspects of the Langlands programme 
started in the 1970's, notably in \cite{harder-75, asai, kurcanov, grunewald-mennicke, cremona-84}. However, up to day, most of the fundamental problems in 
this area are still wide open. In the recent years, the area started receiving attention again and has witnessed a significant amount of activity together with some exciting progress. 

In this article, we aim to discuss several prominent problems of arithmetic nature that surround Bianchi groups together with a survey of existing literature 
including recent work. In agreement with the theme of these proceedings, we will pay attention to computational aspects of these problems, presenting numerical data and concrete examples along the way. We hope that this survey article will contribute to the recently rekindled interest in Bianchi groups.

In Section \ref{basics}, we cover some basic aspects of Bianchi groups that are relevant to our main discussions. In Sections \ref{cohomology} and \ref{hecke}, we discuss the cohomology of Bianchi groups and the Hecke action on the cohomology respectively. In Section \ref{modularforms}, we consider modular forms associated to 
Bianchi groups. In Section \ref{dimension}, we talk about the problem of understanding the size of the cohomology of Bianchi groups. We also discuss the recent asymptotic results of Calegari-Emerton \cite{calegari-emerton}, Finis-Grunewald-Tirao \cite{fgt}, Marshall \cite{marshall} and \c{S}eng\"un-T\"urkelli 
\cite{sengun-turkelli}.  In Section \ref{adelic}, we discuss the adelic locally symmetic spaces that arise in the study of Bianchi groups. 
The subject of Section \ref{elliptic} is the conjectural connection between elliptic curves over imaginary quadratic fields and modular 
forms associated to Bianchi groups. We also talk about the case of higher dimensional abelian varieties and include 
a discussion of elliptic curves with everywhere good reduction. In Section \ref{torsion}, we focus on the torsion in the homology of Bianchi groups. 
We include a discussion of the asymptotic growth of the torsion in light of the recent work of Bergeron and Venkatesh \cite{bergeron-venkatesh}. 
In Sections \ref{galois1} and \ref{galois2}, we consider the (mostly) conjectural connections between 2-dimensional representations of the absolute 
Galois group of imaginary quadratic fields and classes in the cohomology of Bianchi groups. In particular, in the latter section
we make some speculations which relate {\em even} irreducible 2-dimensional representations of the absolute 
Galois group of $\Q$ to the torsion in the cohomology of Bianchi groups.   

There are many significant arithmetic aspects that are left untouched in this survey. 
One of the prominent such aspects is the so called ``$R=T$" problem which has been studied in \cite{bk1, bk2, bk3}.
We also excluded discussions of various $L$-functions and their special values, studied notably in \cite{zhao, urban, ghate, namikawa}, 
and finer arithmetic properties of mod $p$ eigenvalue systems, studied notably in \cite{edixhoven-khare, sengun-turkelli-weight, mohamed-weight}. 

During the final revision phase of this paper, significant developments occured. Calegari and Geraghty made a breakthrough \cite{calegari-geraghty-1, calegari-geraghty-2} 
on the $R=T$ problem using a novel method that allows one go beyond that of Taylor-Wiles. Moreover, Calegari and Venkatesh 
released the first draft of their book project \cite{calegari-venkatesh} on the investigation of integral versions of aspects of the Langlands programme 
carried out in the setting of $\gl$ over imaginary quadratic fields. Unfortunately, due to lack of time, we will not be able to discuss the striking results and 
insights of \cite{calegari-venkatesh}.

\subsection{Notation}
Let $d>0$ be a square-free integer and $K=\Q(\sqrt{-d}) \subset \C$ be an imaginary quadratic field.  Let $D=D_K$ be the discriminant of $K / \Q$. 
We have $D=-4d$ or $D=-d$ if $-d \equiv 2,3$ or $-d \equiv 1$ modulo $4$ respectively. 
The ring of integers $\O = \O_K=\O_d$ of $K$ has a $\Z$-basis $\{ 1, \omega \}$ where 
$\omega = \frac{D+\sqrt{D}}{2}$. The letter $G$ will denote a Bianchi group. The absolute Galois groups of 
$\Q$ and $K$ will be denoted by $G_\Q$ and $G_K$ respectively.

\section{Basic Properties of Bianchi Groups} \label{basics}

\subsection{Geometric Background}
Consider the hyperbolic 3-space $\H$, that is, the unique connected, simply connected Riemannian manifold of dimension 3 with constant sectional curvature $-1$.  A standard model for $\H$ is the upper half space model 
$$\H = \{ (x,y) \in \C \times \R \mid y>0 \}$$
with the metric coming from the line element
$$ds^2 = \dfrac{dx_1^2+dx_2^2+dy^2}{y^2}$$
with $x=x_1+ix_2$. 
Every element $\left ( \begin{smallmatrix} a & b \\ c&d \end{smallmatrix} \right )$ of $\psl(\C)$ acts on $\H$ 
as an orientation-preserving isometry  via the formula 
$$\begin{pmatrix} a & b \\ c&d \end{pmatrix} \cdot (x,y) = 
\Biggl ( \frac{(ax+b)\overline{(cx+d)}+a\bar{c}y^2}{|cx+d|^2+|c|^2y^2},\frac{y}{|cx+d|^2+|c|^2y^2} \Biggr )$$
where $(x,y) \in \H$. It is well-known that every orientation-preserving isometry of $\H$ arises this way. 
The group of {\em all} isometries of $\H$ is generated by $\psl(\C)$ and complex conjugation $(x,y) \mapsto (\bar{x},y)$, which is an orientation-reversing involution. 

As $\O_K$ is a discrete subring of $\C$, the Bianchi group $G$ is a discrete subgroup of the real Lie group $\psl(\C)$. As a result, it acts properly discontinuously on $\H$. The quotients $Y_G =  G \backslash \H$ are hyperbolic 3-orbifolds which are non-compact but of finite volume.  
In fact, we have 
$$vol(Y_{G})=\dfrac{|D_K|^{3/2}}{4\pi^2}\cdot \zeta_K(2),$$
with $\zeta_K$ being the Dedekind zeta function of $K$, as first proved\footnote{Note that a gap in Humbert's proof was filled in by Grunewald and K\"uhnlein in \cite{grunewald-kuhnlein}.} in \cite{humbert} by Humbert in 1919.

It is a classical result that was first observed by Bianchi and first proved by Hurwitz (1892) that the number of cusps 
of $Y_{G}$ is equal to the class number of $K$. The cross-sections of the cusps are 2-tori except when 
$d=1,3$. In these exceptional cases, there are units in $\O_K$ other than $\pm 1$ and as a result the cross-sections 
of the cusps are $2$-folds with $2$-sphere as the underlying manifold.  The orbifold structure of $Y_{G}$ is 
pictured by Hatcher \cite{hatcher} for $D>-100$.

The action of $G$ on $\H$ admits a fundamental domain that is bounded by finitely many geodesic surfaces. 
It follows that $G$ has a finite presentation. 

\subsection{Explicit Fundamental Domains} 
\begin{figure}[htbp]
\begin{center}
\includegraphics[width=0.75\textwidth]{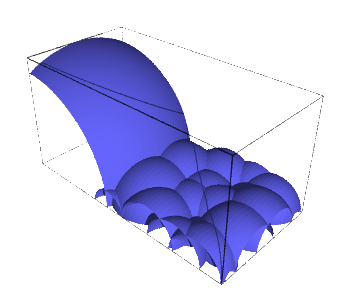}
\end{center}
\caption{The points lying above the blue hemispheres form a fundamental domain for the Bianchi group associated to  $\Q(\sqrt{-39})$.}\label{fig:BD-39}
\end{figure}

\begin{figure}[htbp]
\begin{center}
\includegraphics[width=0.4\textwidth]{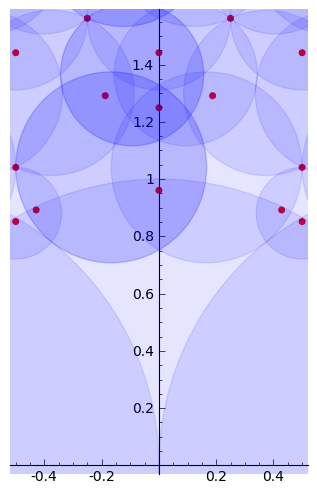}
\end{center}
\caption{Projection on $\C \times \{ 0 \}$ of the {\em floor} of the fundamental domain in Figure~\ref{fig:BD-39}.}\label{fig:Proj-39}
\end{figure}

First explicit fundamental domains were computed by Bianchi \cite{bianchi} in 1890's. In the late 1960's Swan elaborated on Bianchi's work in \cite{swan-68, swan-71}, filling in certain gaps and further developing the theory in a way that is suitable for algorithmic purposes. Detailed accounts of fundamental domains for Bianchi groups 
and Swan's work can be found in 
Section 7 of \cite{egm-book}, Section 6.2.4. of \cite{fine-book} and \cite{bygott, aranes, rahm-thesis}. In the mid 1970's Riley \cite{riley}, building on Swan's work, developed a computer program that produced explicit fundamental domains for Kleinian groups which are geometrically finite (that is, they admit Ford domains with finitely many faces). Detailed discussions of Riley's program can be found in \cite{fine-book, fine-90}. In their recent Ph.D. theses, Aran\'es \cite {aranes} and independently Rahm \cite{rahm-thesis} developed freely accessible computer programs that compute 
fundamental domains for Bianchi groups using Swan's methods. We include here (see Figure~\ref{fig:BD-39}) a plotting of the fundamental domain for the Bianchi group associated to $\Q(\sqrt{-39})$ (of class number 4)  
produced by Aran\'es' Sage \cite{sage} implementation. In Figure~\ref{fig:Proj-39} we see the projection onto the $y=0$ plane of the floor of the fundamental domain, with red points marking the projections
of the singular\footnote{These are the cusps that correspond to non-trivial elements of the class group.} cusps and points in the intersection of 3 or more hemispheres.

Another computer program to this effect was developed by Yasaki in \cite{yasaki} for $\gl(\O_K)$. 
His program uses a different approach that we will sketch later in Section \ref{computingcohomology}. In his master's thesis \cite{page, page-paper}, Page developed a computer program to compute fundamental domains of arithmetically defined Kleinian groups. These groups are defined via quaternion algebras over number fields with a unique complex place (see \cite{maclachlan-reid, egm-book}) and Bianchi groups are among these groups.

\subsection{Explicit Presentations} 
When $O_d$ has Euclidean algorithm (this is the case exactly for $d=1,2,3,7,11$), one can algebraically find a presentation for the associated 
Bianchi group, see work of Cohn \cite{cohn} and Fine \cite{fine-74}. However, for the general case, one needs to use geometry. Once a suitable fundamental domain with an identification of the faces is computed, one can (using, for example, \cite{macbeath}) derive a presentation for the Bianchi group.
Many computations of this sort can be found in the literature, for example, see \cite{swan-71, riley, alperin, schneider, yasaki, page}. By investigating 
the action of Bianchi groups on certain 2-dimensional CW complexes (we will discuss these in Section \ref{computingcohomology}), Fl\"oge computed in \cite{floge} explicit presentations for several Bianchi groups. Many of the presentations that exist in the literature are conveniently compiled in \cite{fgt}.  
A presentation with only two generators for $\psl(\Z[i])$ is given by Brunner in \cite{brunner}. The aforementioned computer program of Page \cite{page-paper} also 
computes presentations for arithmetic Kleinian groups from the geometric data it collects.  

For the sake of illustration, we include here a presentation of the group $\psl(\Z[i])$ that goes back Sansone \cite{sansone}:
$$\left \langle a,b,c,d \mid a^3=b^2=c^3=d^2=(ac)^2=(ad)^2=(bc)^2=(bd)^2=1  \right \rangle.$$

\subsection{Group Theoretic Properties}

In \cite{fine-frohman}, Fine and Frohman proved that, except for $d=1,3$, all Bianchi groups are HNN extensions with amalgamated subgroup $\psl(\Z)$. 
This is achieved by showing that the quotient manifold associated to $\psl(\Z)$ embeds as a non-separating
incompressible surface in $Y_G$. Moreover, they showed that except for $d=3$, all Bianchi groups split as amalgamated free products. 
It is known that for $d=3$, this cannot hold\footnote{In this case, Scarth exhibits in \cite{scarth} a finite-index subgroup which splits as an amalgamated product.}. Every Bianchi group can be described as a graph product of groups, see Brunner et al. \cite{blw}, Mendoza \cite{mendoza} 
and Fl\"oge \cite{floge}. In particular, in the case $d=3$, we get a triangle of groups. A detailed discussion of algebraic aspects of Bianchi groups 
can be found in Fine's book \cite{fine-book}. As an example, we present the following iterated amalgamated free product structure of $\psl(\Z[i])$:
$$\left ( S_3 \star_{_{C_3}} A_4 \right ) \star_{_{(C_2 \star C_3)}} \left ( D_2 \star_{_{C_2}} S_3\right ).$$
Here $S_n,A_n,D_n,C_n,$ represent the symmetric group on $n$ letters, the alternating group on $n$ letters, 
the dihedral group of order $2n$ and the cyclic group of order $n$ respectively. Note that $C_2 \star C_3 \simeq \psl(\Z)$.

In his influential 1970 paper \cite{serre-70}, Serre 
proved that Bianchi groups have finite index subgroups which are not of congruence type\footnote{Recall that a subgroup $\Gamma$ is {\em of congruence type} if it contains the kernel of the surjection $\psl(\O) \rightarrow \psl(\O / I)$ for some ideal $I$ of $\O$. }. Constructions of infinite families of non-congruence subgroups for several Bianchi groups can be found in Britto's paper \cite{britto}. The notion of {\em level} for non-congruence subgroups was discussed by Grunewald and Schwermer in \cite{grunewald-schwermer-99}. A detailed discussion around the Congruence Subgroup Problem for Bianchi groups can be found in \cite{fine-book}, see also \cite{gjz}.

The automorphism groups of Bianchi groups were determined by Smillie and Vogtmann in \cite{smilie-vogtmann}. It turns out that the group 
of outer automorphisms is an elementary abelian 2-group given by complex conjugation, the involution arising from 
$\left (   \begin{smallmatrix} -1 & 0 \\ 0 & 1 \end{smallmatrix} \right )$ and the 2-part of the class group of the underlying imaginary 
quadratic field $K$. Orientation-reversing automorphisms can be used to obtain lower bounds for the size of the cohomology 
of Bianchi groups in certain situations, see \cite{harder-75, rohlfs-85, sengun-turkelli, kionke-schwermer} and Section \ref{dimension}.

\section{Cohomology of Bianchi groups} \label{cohomology}
The cohomology of Bianchi groups, with coefficients in finite dimensional, irreducible, complex representations of the 
ambient real Lie group $\psl(\C)$, is important for number theoretical reasons that we shall discuss later. Let $\Gamma$ 
be a finite index subgroup of a Bianchi group $G$ and let $V$ be a such a representation of $\psl(\C)$ regarded as a representation of $\Gamma$. Then $V$ gives rise to a locally constant sheaf $\mathcal{V}$ on the quotient 3-fold $Y_\Gamma$ associated to $\Gamma$ and it follows from the contractibility of $\H$ that 
$$H^i(\Gamma,V) \simeq H^i(Y_\Gamma,\mathcal{V}).$$
Of course, the same is true for the homology groups. 

A general construction of Borel and Serre \cite{borel-serre} gives\footnote{In fact, the case of symmetric spaces of rank one
was done first by Serre in \cite{serre-70}, to investigate the Congruence Subgroup Problem for Bianchi groups, and served as a prototype for the general construction.} a compact 3-fold $X_\Gamma$ with boundary $\partial X_\Gamma$ such that the interior of $X_\Gamma$ is homeomorphic to $Y_\Gamma$ and the embedding $Y_\Gamma \hookrightarrow X_\Gamma$ is a homotopy equivalence. In particular, we have 
$$H^i(Y_\Gamma,\mathcal{V}) \simeq H^i(X_\Gamma, \overline{\mathcal{V}})$$
where the sheaf $\overline{\mathcal{V}}$ is a suitable extension of $\mathcal{V}$ to $X_\Gamma$. Topologically, $X_\Gamma$ is obtained by attaching at infinity a copy of its cross-section to the cusps of $Y_\Gamma$. Thus except when $d=1,3$, the boundary of $X_\Gamma$ is a finite disjoint union of $2$-tori. For details, we recommend Chapter 2 of the lecture notes of Harder \cite{harder-notes}. 

Consider the natural restriction map 
$$res^i: H^i(X_\Gamma, \overline{\mathcal{V}}) \rightarrow H^i(\partial X_\Gamma, \overline{\mathcal{V}}).$$
It is well known that the rank of the image of $res^1$ is half that of the codomain.  Moreover, the map $res^2$ is onto except 
when $V$ is the trivial representation $\C$, in which case the image has corank $1$. We define the {\bf cuspidal cohomology} 
$H^i_{cusp}(Y_\Gamma,\mathcal{V})$, and thus $H^i_{cusp}(\Gamma,V)$, as the kernel of the restriction map $res^i$. 
By Poincar\'e duality, see \cite{ash-stevens}, we have $H^1_{cusp}(\Gamma,V) \simeq H^2_{cusp}(\Gamma,V^*)$.  

For a ring $R$, let $V_k(R)$ denote the space of homogeneous polynomials of degree $k$ in two variables with coefficients in $R$. 
The space $V_k(R)$  is of rank $k+1$ with basis of monomials $x^k, x^{k-1}y, \hdots, xy^{k-1},y^k$ admitting an action of $M_2(R)$, 
$2\times 2$ matrices with entries in $R$, determined by the rule
$$ \begin{pmatrix} a & b \\ c & d \end{pmatrix} x^{k-m} y^m = (ax+cy)^{k-m} (bx+dy)^m$$
for every $\left ( \begin{smallmatrix} a & b \\ c & d \end{smallmatrix} \right ) \in M_2(R)$. When $k!$ is invertible in $R$, $V_k(R)$ is self-dual.

The finite-dimensional, irreducible, complex representations of $\psl(\C)$ as a {\em real} Lie group are given by 
$$V_{k,\ell}(\C) := V_k(\C) \otimes_{\C} \overline{V_{\ell}(\C)}$$
with $k+\ell$ even, where the overline on the second factor means that the action is twisted by complex conjugation. It follows from the above 
paragraph that these representations are self-dual. Results of Borel and Wallach \cite[Section 3]{borel-wallach} (see also Chapter 3 of \cite{harder-notes}) imply that 
$$H^1_{cusp}(\Gamma,V_{k,\ell}(\C)) = \{ 0 \}$$ 
whenever $k \not= \ell$. 

\subsection{Computing The (Co)Homology} \label{computingcohomology}
For reasons that we will explore later, it is very desirable to explicitly compute the cuspidal (co)homology of Bianchi groups and their congruence subgroups. There are several ways to do this but in essence they all depend on explicit geometric information. We sketch many of them below.

\begin{itemize}
\item The earliest systematic homology computations were carried out by Grunewald \cite{grunewald-mennicke} in 1978. In the core of his approach lies the basic fact that $H_1(\Gamma,\Z) \simeq \Gamma^{ab}$ for any $\Gamma$. Given an explicit presentation of $\Gamma$, one can easily 
compute $\Gamma^{ab}$. Using explicit presentations of Bianchi groups 
$\psl(\O_d)$ for $d=1,2,3$, together with explicit coset representatives, he computed\footnote{Segal records the following in \cite{segal} about these computations of Grunewald: 
\begin{quotation} ``Fritz's gentle manner concealed a steely determination when it came to serious computation. He
was banned for a while from using the Bielefeld University mainframe after his program had
monopolized the entire system: in order to carry out the heavy-duty computation required
for this project, he had devised a routine that managed to bypass the automatic quota
checks."\end{quotation}} explicit presentations (and thus the abelianizations) of congruence subgroups of the form $\Gamma_0(\p)$ for prime ideals $\p$ of residue degree one and norm $\leq 1500$. The original 1978 manuscript is unpublished but some of the results in there were reported in the articles \cite{ghm, egm-80, egm-81}. Recently the author used this technique in \cite{sengun-torsion} and extended the above computations to $d=1,2,3,7,11$ and to prime ideals $\p$ of residue degree one and norm up to 100,000.

\item In 1981 Cremona \cite{cremona-84} adopted the {\em modular symbols} method, developed by Manin (see 
\cite{manin, mazur, cremona-book}) for the classical modular group $\psl(\Z)$, to study the cuspidal part of $H_1(\Gamma,\Q)$ 
for congruence subgroups of the form $\Gamma_0(\a)$ for ideals $\a$ of norm $\leq 500,300,500,200,200$ inside Bianchi groups 
$\psl(\O_d)$ for $d=1,2,3,7,11$ respectively. This technique was adapted by Cremona's students to other Bianchi groups, 
namely, by Whitley \cite{whitley} to the remaining class number one cases $d=19,43,67,163$, 
by Bygott \cite{bygott} to $d=5$ (class number 2) and by Lingham \cite{lingham} to $d=23,31$ (both class number 3). 

Let $Y^*_\Gamma=\Gamma \backslash \H^*$ be the Satake compactification of $Y_\Gamma = \Gamma \backslash \H$ 
where $\H^*=\H \cup \mathbb{P}^1(K)$. Then the homology $H_1(Y^*_\Gamma, \Q)$ is isomorphic to the cuspidal part of 
$H_1(Y_\Gamma,\Q) \simeq H_1(\Gamma, \Q)$. It can be shown that $H_1(Y^*_\Gamma, \Q)$ can be generated by translates 
of finitely many geodesics (which are explicitly determinable) with both ends in $\mathbb{P}^1(K)$ by coset representatives of $\Gamma$ in $G$. 
One needs explicit geometric information to compute the complete system of relations between these generators. For details see \cite{cremona-book, 
cremona-notes}.

\item  
Given any torsion-free finite index subgroup $\Gamma$ of a Bianchi group $G$, we have $H^i(\Gamma, M)=0$ for every 
finite dimensional $\Gamma$-module $M$ and $i>2$, in other words, the virtual cohomological dimension of $G$ is 2. 
Thus from a cohomological perspective, our geometric playground $\H$ is not of optimal dimension. 
Using reduction theory one can construct a 2-dimensional CW-complex $\mathcal{C}$ on which $G$ acts cellularly such that 
$G \backslash \mathcal{C}$ is finite. Given such a $CW$-complex, one gets explicit descriptions of all (co)homology groups 
in degrees $0,1,2$ by feeding the stabilizers of the cells and the incidence relations between the cells into the so called {\em equivariant 
(co)homology spectral sequence}. 

Such CW complexes were first constructed by Mendoza \cite{mendoza} and Fl\"oge \cite{floge} as $G$-invariant deformation retracts in $\H$, also known as {\em spines}. Mendoza considered the groups $\psl(\O_d)$, $\pgl(\O_d)$ for $d=1,2,3,7,11$ while Fl\"oge worked only with $\psl(\O_d)$ 
but for more cases, namely for $d=1,2,3,7,11,5,6,10$. Mendoza retracts away from the cusps whereas Fl\"oge adjoins the singular cusps to his 
CW-complex. In particular, the stabilizers of the 0-cells given by these singular cusps have infinite stabilizer groups. In such cases, 
the spectral sequence degenerates on the $E^3$-page.

These CW-complexes were used by Schwermer and Vogtmann \cite{schwermer-vogtmann} to compute the integral homology of $\psl(\O_d)$ 
for $d=1,2,3,7,11$ in all\footnote{Note that beyond degree two, the homology groups are exclusively torsion at the primes 2 and 3. Their explicit 
structure depends on the number of conjugacy classes of the finite subgroups of the Bianchi group, see Rahm \cite{rahm-farrell}.} degrees. The ring structure of the integral homology was computed by Berkove \cite{berkove} using a similar approach. Later Vogtmann \cite{vogtmann} wrote a computer program to compute the spines using Mendoza's method and computed the rational second homology of $\psl(\O_d)$ with discriminant $D>-100$. A similar computer program was written by Scheutzow \cite{scheutzow} who extended the computations of Vogtmann to both $\psl(\O_d)$ and $\pgl(\O_d)$ for all $K$ with discriminant $D \geq -260$ and also for $d=643$. The author used results of Mendoza for extensive computations in \cite{sengun-torsion}. Rahm in his 2010 thesis \cite{rahm-thesis} wrote a computer program to compute spines using the approach of Fl\"oge. This program was used by Rahm and Fuchs in \cite{rahm-fuchs} to compute the integral homology in all degrees of $\psl(\O_d)$ for  $d=5,6,10,13,15$, all of class number 2. Recently Rahm and the author used this program in \cite{rahm-sengun} to compute more than 4900 instances of $H^2(\psl(\O_d),V_{k,k}(\C))$ for more than 180 values of $d$. 
 
Finally, let us note that in \cite{alperin}, Alperin constructed such a CW complex for $\psl(\O_3)$ using a method that is different than 
those of Mendoza and Fl\"oge.

\item Given a finitely presented group $G$, one can compute its first homology and first cohomology using the so called {\em Fox Calculus}, see \cite{fox}. The rough idea is that a 1-cocycle $f:G \rightarrow M$ can be determined uniquely by its values on the generators. It follows that 
the set of 1-cocycles can be identified with elements of $M^n$ (here $n$ is the number of generators) satisfying the system of linear 
equations that arise from the relators of $G$. This approach has been used to compute the first cohomology of several Bianchi groups in 
\cite{priplata, sengun-thesis, fgt}.

\item Recently Yasaki \cite{yasaki} developed, building on work of Gunnells \cite{gunnells}, a computer program which computes 
$H^1_{cusp}(\Gamma_0(\mathfrak{n}),\C)$ with $\Gamma_0(\mathfrak{n}) \subset \gl(\O_K)$ that works with any given imaginary 
quadratic field $K$. This program is now part of Magma \cite{magma}. The space of positive definite binary Hermitian forms over $K$ forms an open cone in a real vector space. There is a natural decomposition of this cone into polyhedral cones corresponding to the facets of the Voronoi polyhedron, 
see \cite{gunnells}. The top-dimensional polyhedral cones of the decomposition correspond to perfect forms and descend to ideal polytopes in $\H$. 
In \cite{yasaki} Yasaki gives a method to compute an initial perfect form, which is the necessary input for a general algorithm by 
Gunnells \cite{gunnells} that classifies perfect forms under the action of $\gl(\O_K)$. Note that this is a special case of a 
powerful method that is able to capture the cohomology of $\textrm{GL}_n$ over number fields, see \cite{agm,ghy,gy}. For 
a very nice exposition, see Gunnells' lecture notes \cite{gun13} in these proceedings. 

\item In her recent paper \cite{torrey}, Torrey developed a computer program to compute with the cohomology of $\gl(\Z[i])$ and its congruence subgroups of type $\Gamma_0$ using the {\em Borel-Serre duality}, which provides a Hecke module isomorphism $H^2(\Gamma,V) \simeq H_0(\Gamma, St \otimes V)$. Here $St$ denotes the {\em Steinberg module} $H^2(\Gamma, \C[\Gamma])$, which can be described in terms of modular symbols via Ash's work \cite{ash-94}. 

\end{itemize}


\section{Hecke Action} \label{hecke}

The cohomology groups $H^i(\Gamma,V) \simeq H^i(\Gamma \backslash \H, \mathcal{V})$ come equipped with an infinite commuting family
 of endomorphisms, called {\bf Hecke operators}, which are crucial for connections with number theory. These come from correspondences 
which arise from elements of the commensurators of Bianchi groups. 

Let $\Gamma$ be a finite index subgroup of a Bianchi group $\psl(\O_K)$. Given $g \in \pgl(K)$, let us consider 3-folds 
$M,M_g,M^g$ associated to $\Gamma, \Gamma_g := \Gamma \cap g \Gamma g^{-1}, 
\Gamma^g :=\Gamma \cap g^{-1} \Gamma g $ 
respectively. We have {\em finite} coverings, induced by inclusion of fundamental groups, 
$$r_g:M_g \rightarrow M, \ \ \ \ \ r^g: M^g \rightarrow M$$
and an isometry 
$$\tau : M_g \rightarrow M^g$$
induced by the isomorphism between $\Gamma_g$ and $\Gamma^g$ given by conjugation by $g$.
The composition $s_g := r^g \circ \tau$ gives us a second finite covering from $M_g$ to $M$. 
The coverings $r_g$ induce linear maps between the homology groups 
$$r_g^* : H^i(M,\mathcal{V}) \rightarrow H^i(M_g,\mathcal{V}).$$
The process of summing the finitely many preimages in $M_g$ of a point of $M$ under $s_g$ leads to 
$$s_g^* : H^i(M_g,\mathcal{V}) \rightarrow H^i(M,\mathcal{V}).$$
Note that $s_g^*$ is equivalent to the composition
$$H^i(M_g,\mathcal{V}) \rightarrow H^i(M^g,\mathcal{V}) \rightarrow H^i(M,\mathcal{V})$$
where the first arrow is induced by $\tau$, 
and the second arrow is simply the corestriction map (which corresponds to the ``transfer map" of group cohomology). We define the {\bf Hecke operator} $T_g$ associated to $g \in \pgl(K)$ as the composition 
$$T_g:= s_g^* \circ r_g^* : H^i(M,\mathcal{V}) \rightarrow H^i(M,\mathcal{V}).$$  
There is a notion of isomorphism of Hecke operators that we shall not present. It turns out that up to isomorphism, 
$T_g$ depends only on the double coset $\Gamma g \Gamma$.

One can define Hecke operators using the above process for homology groups as well. 
It might provide insight to look at the situation from the perspective of $H_2(M,\Z)$ (which is isomorphic to compactly supported cohomology $H^1_c(M,\Z)$ by duality, see \cite{ash-stevens}). If an embedded surface $S \subset M$ represents $\omega \in H_2(M,\Z)$, then the immersed surface $r_g(s_g^{-1}(S))$ represents $T_g\omega \in H_2(M,\Z)$. 

Given an element $\pi \in \O_K$, the {\em Hecke operators $T_\pi,S_\pi$ associated to $\pi$} are defined as the Hecke operators associated to the matrices 
$( \begin{smallmatrix}  \pi & 0 \\ 0 &1 \end{smallmatrix} )$ and $( \begin{smallmatrix}  \pi & 0 \\ 0 &\pi \end{smallmatrix} )$ in $\gl(K)$.
Let $\mathbb{T}$ denote the subalgebra of the endomorphism algebra of $H^i(\Gamma,V)$ that is generated by the Hecke operators $T_{\pi},S_\pi$ where 
$\pi$ is a prime element of $\O_K$ such that $(\pi)$ is coprime to the level of $\Gamma$. 

It is well-known that the Hecke algebra $\mathbb{T}$ is commutative, semi-simple and stabilizes the cuspidal part $H^i_{cusp}(\Gamma,V)$. The duality between $H^1_{cusp}(\Gamma,V)$ and $H^2_{cusp}(\Gamma,V)$ is $\mathbb{T}$-equivariant, see \cite{ash-stevens}. So for number theoretic purposes, 
they carry the same ``information", that is, eigenvalue systems. An {\bf eigenvalue system} of $\mathbb{T}$ with values in a ring $R$ is a ring homomorphism $\Phi : \mathbb{T} \rightarrow R$. We say that an eigenvalue system $\Phi$ {\em occurs} in the $R\mathbb{T}$-module $A$ if there is a nonzero element $a \in A$ such that $Ta=\Phi(T)a$ for all $T$ in $\mathbb{T}$. For an eigenvalue system $\Phi$ that occurs in $H^i(\Gamma, V_{k,\ell}(\C))$, the image $\Phi(\mathbb{T})$ is not just a collection complex numbers. It follows from the isomorphism 
$H^i(\Gamma, V_{k,\ell}(\C)) \simeq H^i(\Gamma,V_{k,\ell}(K)) \otimes_K \C$ that $\Phi(\mathbb{T})$ lies in an algebraic number field $L$. In fact, after 
a simple normalization, $\Phi(\mathbb{T})$ lies in the ring of integers $\O_L$ of $L$. 

These eigenvalue systems are closely connected to 2-dimensional representations of the absolute Galois group of $K$ and to abelian varieties 
as we shall discuss later on. However when considering such connections, it is necessary that the Hecke operators are indexed by the {\em ideals} of $\O_K$. This amounts to 
having $T_\pi = T_{-\pi}$ for every element $\pi \in \O_K$. In order to achieve this, one needs to work with the cohomology of $\pgl(\O_K)$.

Finally let us remark that the action of Hecke operators in the case when $\Gamma$ is not a congruence subgroup has been studied by Zarghani in \cite{zarghani}.

\subsection{Computing the Hecke Action}
 
Computing the Hecke action on the (co)homology of arithmetic groups is a highly difficult task in general (see \cite{gunnells-hecke-1, 
gunnells-hecke-2, gunnells-mcconnell}). In the case of linear groups over 
a number field with a Euclidean algorithm, computing the Hecke action on the top degree cohomology boils down to a continued 
fraction algorithm (see \cite{ash-lee}). 

In the case of $\psl(\Z)$, Manin shows in \cite{manin} that the first rational homology can be generated by certain distinguished 
modular symbols and the Hecke action on these distinguished symbols can be computed using the continued fractions algorithm. This approach 
can be adapted to the case of Bianchi groups over the Euclidean imaginary quadratic fields $K$ as was done by Cremona in \cite{cremona-84}.  In \cite{whitley}, Whitley developed a geometric counter-part, called {\em the pseudo-Euclidean algorithm}, of the continued fractions algorithm that also works in the cases where $K$ no longer has Euclidean
 algorithm. The algorithm finds a matrix that sends a given cusp either to the cusp $\infty$ or to the ``singular" point in its orbit 
using the geometry of the Bianchi group (a suitable fundamental domain with the list of hemi-spheres 
covering its floor should be pre-computed (see \cite[p.87]{aranes}). The method was further worked on by Bygott in \cite{bygott} and used 
in \cite{bygott, lingham}. 

In practice, when one works with modular symbols, it is more convenient to employ the so-called {\em Manin symbols} (see \cite[Section 2.2]{cremona-book}). In the case of $\psl(\Z)$, the action of Hecke operators on Manin symbols was described by Merel \cite{merel}, see also Zagier \cite{zagier}. In his 2011 Ph.D. thesis \cite{mohamed}, Mohamed described the action of Hecke operators on Manin symbols in the setting of Bianchi 
groups associated to imaginary quadratic fields of class number one. Mohamed's description allows one to compute 
Hecke operators on the second cohomology of Bianchi groups as well since the previously discussed method of computing the second cohomology 
via spines boils down to computing certain co-invariants, which can be interpreted as Manin symbols, see \cite{mohamed, mohamed-12}. Indeed, this was the 
method that we used to compute with Hecke operators in the numerical examples of Section \ref{even}.

Let us finally note that the above mentioned process of using continued fractions algorithm to compute Hecke operators is equivalent to writing 
a given element as a word in a certain set of generators of the group in the cases of $\psl(\Z)$ and the Bianchi groups associated to the Euclidean imaginary quadratic fields, see \cite[Section 9]{rademacher} and \cite{series}. 


\section{Bianchi Modular Forms and Cohomology} \label{modularforms}

Just as the cohomology of the classical modular group $\psl(\Z)$ is central to 
the theory of elliptic modular forms, the cohomology of Bianchi groups is central to the study of 
Bianchi modular forms, that is, modular forms over imaginary quadratic fields. These are {\em vector valued} real analytic functions on $\H$ with certain transformation properties satisfying certain differential equations and growth properties. 
Detailed discussions can be found in Kubota \cite{kubota-64}, Friedberg \cite{friedberg-1, friedberg-notes} and Bygott \cite{bygott}. 
Let us give here a quick definition for the case of weight 2 from the perspective of vector valued harmonic differential 1-forms. Note that in the cases of weights 0 and 1, the automorphic forms are not {\em cohomological}, that is, 
are not captured in the cohomology of Bianchi groups.

Given $\gamma=(\begin{smallmatrix} a & b \\ c&d \end{smallmatrix}) \in \psl(\C)$ and $z=(x,y) \in \H$, 
let us introduce the multiplier system
$$J(\gamma, z):= \begin{pmatrix} cx+d & -cy \\ \bar{c}y & \overline{cx+d}\end{pmatrix}$$

Given a function $F: \H \rightarrow \C^{k+1}$ and $\gamma \in \psl(\C)$, we define the {\em slash operator} 
$$(F |_k\gamma)(z):=\sigma^k(J(\gamma, z)^{-1}) F(\gamma z)$$
where $\sigma^k$ is the {\em symmetric $k^{th}$ power} of the standard representation of $\psl(\C)$ on $\C^2$.

In the case $k=2$, we have $F: \H \rightarrow \C^3$ and
$$(F |_k\gamma)(z)=\dfrac{1}{|r|^2+|s|^2} 
\begin{pmatrix}  \bar{r}^2 & 2\bar{r}s & s^2 \\ 
                          -\bar{r}\bar{s} & |r|^2-|s|^2 & rs \\ 
                          \bar{s}^2 & -2r\bar{s} & r^2 \end{pmatrix} F(\gamma z)$$
where $\gamma=( \begin{smallmatrix} a &b \\ c&d \end{smallmatrix})$ and $r=cx+d$ and $s=cy$.

The 1-forms $\beta_1:=-\frac{dx}{y}, \beta_2:= \frac{dy}{y}, \beta_3:=\frac{d\bar{x}}{y} $ form a basis of differential 1-forms on $\H$. 
A differential form $\omega$ is {\em harmonic} if $\Delta \omega =0$ where $\Delta=d \circ \delta + \delta \circ d$ 
is the usual Laplacian with $d$ being the exterior derivative and $\delta$ the codifferential operator. 
Then $\psl(\C)$ acts on the space of differential 1-forms as 
$$\gamma \cdot {}^t(\beta_1,\beta_2,\beta_3)_{(z)} = \sigma^2(J(\gamma,z)){}^t(\beta_1,\beta_2,\beta_3)_{(z)} .$$
A {\bf weight $2$ cuspidal Bianchi modular form} for a congruence subgroup $\Gamma \subset \psl(\O_K)$ is a real analytic function $F=(F_1,F_2,F_3) : \H \rightarrow \C^3$ 
with the following properties.
\begin{enumerate}
\item $F_1\beta_1 + F_2 \beta_2+F_3\beta_3$ is a harmonic differential 1-form on $\H$ that is $\Gamma$-invariant,
\item $\int_{\C / \O_K} (F|\gamma )(x,y) dx = 0$ for every $\gamma \in \psl(\O_K).$
\end{enumerate}
The last condition is equivalent to saying that the constant coefficient in the {\em Fourier-Bessel expansion} of $F|\gamma$ is equal to zero for every $\gamma \in \psl(\O_K)$. Let us explain this. As $\Gamma$ is non-cocompact, it contains parabolic elements. The $\Gamma$-invariance of $F$, which is implicit in condition (1), implies that $F$ is a periodic function in the $x=(x_1,x_2)$-variable. It follows that $F$ has a Fourier-Bessel expansion of the form 
$$F(x,y)=\sum_{0 \not =\alpha \in \O_K}c(\alpha) y^2 \mathbb{K}\left ( \dfrac{4\pi|\alpha|y}{\sqrt{|\triangle|}} \right )
                        \psi\left (\dfrac{\alpha x}{\sqrt{|\triangle|}} \right )$$
where 
$$\psi(x)=e^{2\pi(x+\bar{x})}$$
and  
$$\mathbb{K}(y)=\left ( -\dfrac{i}{2}K_1(y),K_0(y),\dfrac{i}{2}K_1(y) \right)$$
with $K_0,K_1$ are the hyperbolic Bessel functions satisfying the differential equation 
$$\dfrac{dK_j}{dy^2}+\dfrac{1}{y}\dfrac{dK_j}{dy}-\left ( 1+\dfrac{1}{y^{2j}}\right )K_j = 0, \ \ \ \ j=0,1$$
and decreases rapidly at infinity. 

The space of weight k cuspidal modular forms for a fixed $\Gamma$ is a finite dimensional complex vector space which we will denote with $S_k(\Gamma)$. This space comes equipped with an infinite commuting family of endomorphisms called {\em Hecke operators}, indexed by 
elements of $\O_K$. The {\bf generalized Eichler-Shimura Isomorphism}\footnote{This is also known as the {\em Eichler-Shimura-Harder Isomorphism}. However Prof.Harder advised 
us not to use this name.}, as studied by Harder in \cite{harder-75, harder-87}, provides us with an explicit 
Hecke module isomorphism
$$H^1_{cusp}(\Gamma, V_{k,k}(\C)) \simeq S_{k+2}(\Gamma).$$  
See also Kur\v{c}anov \cite{kurcanov} for the case of $k=0$.


\section{The Dimension Problem}  \label{dimension}
Understanding the behaviour of the dimension of cohomology of congruence subgroups of Bianchi groups is an old 
open problem. Up to now, we do not have any sort of closed formula for the dimension.

It was observed by Swan \cite{swan-68, swan-71} for several values of $d$ that the dimension of 
$H^1_{cusp}(\psl(\O_d),\C) $ is $0$. Mennicke showed 
(see the end of \cite{swan-71}) that this was not true in general; he computed that the dimension of 
$H^1_{cusp}(\psl(\O_{10}),\C)$ is $1$. 

Zimmert in his 1971 Diplomarbeit \cite{zimmert} showed, via topological means, that there is a set $Z(K)$, today known as the {\em Zimmert Set}, 
associated to every imaginary quadratic field $K$ such that 
$$\# Z(K) \leq \dim H^1_{cusp}(\psl(\O_K), \C).$$
Grunewald and Schwermer \cite{grunewald-schwermer-81-2} modified the purely topological method of Zimmert and proved that 
$$\dim H^1_{cusp}(\psl(\O_K),\C) \rightarrow \infty$$
as the discriminant of $K$ grows. Using the same method, they also proved in 
\cite{grunewald-schwermer-81-3} that every finite index subgroup of a Bianchi group $\psl(\O)$ has finite index (torsion-free) 
subgroups with arbitrarily large cuspidal cohomology.

In 1984 Rohlfs showed in \cite{rohlfs-85}, using an idea that goes back to Harder \cite{harder-75}, that 
$$\dfrac{1}{24} \varphi(d) - \dfrac{1}{4} - \dfrac{h(d)}{2} \leq \ \dim H^1_{cusp}(\psl(\O_d),\C) $$
where $\varphi(d)$ is the Euler $\varphi$-function and $h(d)$ is the class number of $\O_d$ (see \cite{schwermer} for a nice 
exposition of Rohlfs' method). Using this estimate together with the results of Grunewald and Schwermer, he gave a finite set $S$ such that if 
the dimension of $H^1_{cusp}(\psl(\O_d),\C)$ is $0$ then $d \in S$. Kr\"amer \cite{kraemer} and independently Vogtmann \cite{vogtmann} proved later 
that   $\dim H^1_{cusp}(Y,\C) = 0$ only for   
$$d \in \{ 1, 2, 3, 5, 6, 7, 11, 15, 19, 23, 31, 39, 47, 71 \}.$$ 

In fact Kr\"amer produced a lower bound for $H^1_{cusp}(\psl(\O_d),\C)$ that is 
sharper than that of Rohlfs. Methods of Rohlfs and Kr\"amer were developed further by Blume-Nienhaus in \cite{blume} 
where he computed lower bounds $H^1_{cusp}(\Gamma,\C)$ for  groups $\Gamma$ of elements of norm one in 
maximal orders in $M_2(K)$ (see Section \ref{adelic} for these orders).

When we consider congruence subgroups of the form $\Gamma_0(\p)$ with prime $\p$ of residue degree one, 
we see an interesting picture. Extensive data collected by the author in \cite{sengun-torsion} show that roughly speaking 
$90 \%$ of the time, $H^1_{cusp}(\Gamma_0(\p),\C)$ vanishes. However, the ideals $\p$ for which the cuspidal cohomology 
does not vanish seem to have a uniform distribution. We make the following conjecture which is, as we shall discuss in Section 
\ref{elliptic}, is related to the existence of abelian varieties over imaginary quadratic fields with special types of 
endomorphism rings and restricted ramification properties.

\begin{svgraybox}
\begin{conjecture} For every $K$, there are infinitely many prime ideals $\p$ of residue degree one in $\O_K$ such that 
 $$H^1_{cusp}(\Gamma_0(\p),\C) \not= 0.$$
\end{conjecture}
\end{svgraybox}  

At this point a natural question arises.

\begin{svgraybox}
\begin{question} Is there any criterion that allows one to pre-determine for a given prime ideal $\p$ of residue degree one,  whether 
 $$H^1_{cusp}(\Gamma_0(\p),\C) \not= 0 \ \ ?$$
\end{question}
\end{svgraybox}  

Let us now fix the level and consider higher weights. Let $H^1_{bc}(\psl(\O_K),V_{k,k}(\C))$ denote the subspace of 
$H^1_{cusp}(\psl(\O_K),V_{k,k}(\C))$ that consists of classes which correspond to cuspidal Bianchi modular forms 
arising from elliptic modular forms via base change (and their twists) or correspond to cuspidal Bianchi modular forms arising from quadratic characters 
associated with certain quadratic extensions $L/K$ (see \cite{fgt} for details). 
Numerical data, collected in \cite{calegari-mazur, fgt} and remarkably extended in \cite{rahm-sengun}, 
show that {\em generically} this distinguished subspace exhausts all of the cuspidal cohomology.

Accordingly, we include here the next question which was posed in \cite{fgt}.  
\begin{svgraybox}
\begin{question} \label{bc} For a fixed $K$, are there only finitely many weights $k$ such that 
 $$H^1_{bc}(\psl(\O_K),V_{k,k}(\C)) \not= H^1_{cusp}(\psl(\O_K),V_{k,k}(\C)) \ \ ?$$
\end{question}
\end{svgraybox}  

Note that a formula for the dimension of $H^1_{bc}(\psl(\O_K),V_{k,k}(\C))$ is computed in \cite{fgt}. 

\subsection{Asymptotic Behaviour}
Recently there has been progress in understanding the asymptotic behaviour of the dimension.
Let us start with results in the ``vertical" direction. Calegari and Emerton considered in \cite{calegari-emerton} 
how the size of the cohomology, with a fixed coefficient module, varied in a tower of arithmetic groups. 
Their general result when applied to our situation gives the following. Let $\Gamma(\p^n)$ denote the 
principal congruence subgroup of level $\p^n$ of a Bianchi group $\psl(\O)$ where $\p$ is an unramified prime ideal 
of $\O$. Fix $V(\C)$. If the residue degree of $\p$ is one, then 
$$ \dim H^1(\Gamma(\p^n),V(\C)) \ll p^{2n},$$
and if the residue degree of $\p$ is two, then 
$$ \dim H^1(\Gamma(\p^n),V(\C)) \ll p^{5n}$$
as $n$ increases. Note that the trivial upper bounds, which are given by the index of $\Gamma(\p^n)$ in 
the ambient Bianchi group, are $p^{3n}$ and $p^{6n}$ respectively.  It is natural 
to look at these asymptotics from the perspective of the volume which is a topological 
invariant in our setting. Observe that the volume of $Y_{\Gamma(\p^n)}$ is given by 
a constant times the index of $\Gamma(\p^n)$ in the Bianchi group. Thus 
asymptotically, the trivial asymptotic upper bound for the above cohomology groups 
is {\em linear} in the volume and the above upper bounds of Calegari and Emerton can be interpreted 
as {\em sublinear}.

For the ``horizontal direction", let $\Gamma$ be a congruence subgroup of a Bianchi group. 
Finis, Grunewald and Tirao showed in \cite{fgt} that
$$ k \ll \dim H^1_{cusp}(\Gamma, V_{k,k}(\C)) \ll \frac{k^2}{\log k}$$
as $k$ increases. The inequality on the left is proven only for the case $\Gamma=\psl(\O)$. We note that the trivial asymptotic upper bound, 
which is given by the dimension of $V_{k,k}(\C)$, is $k^2$. A recent result of Marshall in \cite{marshall}, which uses ideas of Calegari and Emerton mentioned 
above, when applied to our situation, improves the upper bound by a power:
$$ \dim H^1_{cusp}(\Gamma, V_{k,k}(\C)) \ll_{\epsilon} \ k^{5/3+\epsilon}$$
as $k$ increases. 

Building on \cite{rohlfs-78} and \cite{blume}, T\"urkelli and the author derived, in \cite{sengun-turkelli}, explicit lower bounds which give the following asymptotics. Let $p$ be a rational prime that is unramified in $K$ and let $\Gamma(p^n)$ denote the  principal congruence subgroup of level $(p^n)$ of a Bianchi group $\psl(\O)$. Then 
$$  k \ll \dim H^1_{cusp}(\Gamma(p^n),V_{k,k}(\C))$$
as $k$ increases and $n$ is fixed. Assume further that $K$ is of class number one and that $p$ is inert in $K$. Then 
$$p^{3n} \ll \dim H^1_{cusp} (\Gamma(p^n),\C)$$ 
as $n$ increases.


\section{The Adelic Approach}  \label{adelic}

We have seen that there is an infinite collection of commuting linear operators, called Hecke operators, acting on the cohomology 
of Bianchi groups. In particular, when we work with $\pgl(\O_K)$, the Hecke operators are associated to principal ideals of $\O_K$. From 
an arithmetic point of view, this does not fully reflect the complexity of the base field $K$ as the principal ideals only correpsond to 
the identity element of the class group $\mathcal{C}_K$ of $K$. In this section we discuss the adelic approach for $\gl$
which succeeds in assigning a Hecke operator to {\em every} ideal of $\O_K$.

To this end, for a given non-zero ideal $\b$ of $\O_K$, let us introduce 
$$M(\O_K,\b) :=\left \{ \left ( \begin{smallmatrix} a & b \\ c&d \end{smallmatrix} \right )  \in M_2(K) \mid a,d \in \O_K,\  b \in \b^{-1}, \ c \in \b \right \}$$
where $M_2(K)$ is the quaternion algebra of $2 \times 2$ matrices with entries in $K$. It is well-known 
that $M(\O_K,\b)$ are maximal orders in $M_2(K)$ and every maximal order of $M_2(K)$ is conjugate to 
an $M(\O_K,\b)$ by some element in $\gl(K)$. Moreover, the orders $M(\O_K,\a)$ and $M(\O_K,\b)$ 
are conjugate to each other if and only if $\a$ and $\b$ represent the same class in $\mathcal{C}_K / \mathcal{C}_K^2$.

The groups 
$$\gl(O_K,\b) := \textrm{GL}_1(M(\O_K,\b)) $$
are arithmetic subgroups of $\gl(\C)$, commensurable with $\gl(\O_K,\O_K)$, which is simply $\gl(\O_K)$.

Let $\mathbb{A},\mathbb{A}_f$ denote the rings of adeles and of finite adeles of $K$ and  
$\widehat{\O_K}$ denote the profinite completion of $O_K$. For an open subgroup $\mathcal{K}$ of $\gl(\widehat{\O_K})$, we form the adelic locally symmetric space
$$Y_{\mathbb{A}}(\mathcal{K}):= \gl(K) \backslash \left ( \gl(\mathbb{A}_f) / \mathcal{K} \right ) \times \H .$$
It follows from the strong approximation theorem that the adelic space $Y_{\mathbb{A}}(\mathcal{K})$ is the disjoint union of 
arithmetic hyperbolic 3-folds. More precisely, assume that the determinant map $\mathcal{K} \rightarrow \widehat{O_K}^*$ is surjective. Then 
we can choose ideals $\b_1, \hdots, \b_h$ representing the elements of $\mathcal{C}_K$ with $\b_1=\O_K$ such that
$$ Y_{\mathbb{A}}(\mathcal{K}) = \bigsqcup \Gamma_j \backslash \H$$
where $\Gamma_j = \gl(\O_K, \b_j) \cap \mathcal{K}$ for $j=1,\hdots ,h$. It is interesting to note that the arithmetic groups $\Gamma_j$ are all congruence 
subgroups. The adelic construction does not ``see" the non-congruence arithmetic subgroups of $\gl(K)$, see, for example, \cite[Prop.4.1]{milne}.

\subsection{Adelic Hecke Action}
Let us now see how we can construct a Hecke operator for every ideal of $\O_K$. For our purposes, it suffices to 
work with the level group $\mathcal{K}=\mathcal{K}_1(\mathfrak{n})$ given by
$$\left \{ \left ( \begin{smallmatrix} a & b \\ c&d \end{smallmatrix} \right )  \in \gl(\widehat{\O_K}) \mid c,d-1 \in \mathfrak{n}\widehat{\O_K} \right \}$$
for some fixed ideal $\mathfrak{n}$ of $\O_K$. Notice that $\gl(\O_K) \cap \mathcal{K}_1(\mathfrak{n})$ is the standard 
congruence subgroup $\Gamma_1(\mathfrak{n})$. Now fix a prime ideal $\mathfrak{p}$ that does not divide $\mathfrak{n}$. 
Let $\varpi_{\p}$ denote a uniformizer of the completion $K_{\p}$ of $K$ at $\p$. Then the matrices 
$\left ( \begin{smallmatrix} \varpi_{\p} & 0 \\ 0&1 \end{smallmatrix} \right ), 
\left ( \begin{smallmatrix} \varpi_{\p} & 0 \\ 0&\varpi_{\p} \end{smallmatrix} \right )$ give rise to 
{\em Hecke operators $T_{\p},S_\p$ associated to the prime ideal $\p$} acting 
on $H^i(Y_{\mathbb{A}}(\mathcal{K}), \mathcal{V})$ as was done in Section \ref{hecke} (see \cite{taylor-95} for details). 
These Hecke operators satisfy nice multiplicative properties which allow one to associate a Hecke operator to any ideal that is coprime $\mathfrak{n}$. 
The Hecke algebra $\mathbb{T}_1(\mathfrak{n})$ and the eigenvalue systems on it are defined just as in Section \ref{hecke}.

It follows from the above that the cohomology of the adelic space decomposes as
$$H^i(Y_{\mathbb{A}}(\mathcal{K}), \mathcal{V}) = \bigoplus_{j=1}^h H^i(\Gamma_j \backslash \H, \mathcal{V}).$$
A closer inspection shows that Hecke operators $T_\p$ do not, in general, stabilize the cohomology groups of the connected components 
$\Gamma_j \backslash \H$ of $Y_{\mathbb{A}}(\mathcal{K})$. In fact, $T_\p$ stabilizes these cohomology groups if and only if 
$\p$ is a principal ideal. This is why we were able to define Hecke operators associated to only principal ideals on the cohomology 
of Bianchi groups in Section \ref{hecke}.  Nevertheless, an eigenvalue system on $\mathbb{T}_1(\mathfrak{n})$ can be completely 
recovered from its restriction to the subalgebra $\mathbb{T}$ generated by the Hecke operators associated to principal prime ideals of $\O_K$, 
see \cite{bygott, lingham, cremona-pre} and also \cite{styer}). In fact
two eigenvalue systems have the same restriction to the subalgebra $\mathbb{T}$ if and only if 
one is the twist of the other by a character of $\mathcal{C}_K$. This fact is very convenient for computational purposes as it 
allows us to just work with the cohomology of the congruence subgroup $\Gamma_1(\mathfrak{n})$ of $\gl(\O_K)$ 
even when $K$ has non-trivial class group.

\subsection{Adelic Bianchi Modular Forms}
A cuspidal {\em adelic Bianchi modular form} $F$ of level $\mathcal{K}$ and weight $k+2$ is a vector-valued function on $\gl(\mathbb{A})$ with several 
special properties that we shall not spell out here. However, in accordance with the decomposition of the adelic space 
into connected components, $F$ can be identified with an $h$-tuple 
$(f_1,\hdots,f_h)$ where each $f_j$ is a vector valued real analytic function defined on $\H$ with nice transformation 
properties, dictated by the weight, under the action of $\Gamma_j$ and with controlled growth towards the cusps. In particular, 
$f_1$ is a cuspidal Bianchi modular form of level $\Gamma_1$ and weight $k+2$ as we described in Section \ref{modularforms}. 
Detailed discussions of adelic Bianchi modular forms and their identifications with (tuples of) non-adelic ones can be found in 
Bygott \cite{bygott} and Ghate \cite{ghate}. Adelic Bianchi modular forms can also be studied from the perspective of automorphic representations. 
Detailed discussions within this setting can be found in Harder \cite{harder-87}, Urban \cite{urban} and Berger \cite{berger-thesis, berger-08}.

{\bf From now on}, when we talk about a Bianchi modular form, we shall mean an adelic one for $\gl(\mathbb{A})$, unless stated otherwise.

\section{Connections with Galois Representations I} \label{galois1}

Results of Harris-Soudry-Taylor \cite{hst}, Taylor \cite{taylor-94} and Berger-Harcos \cite{berger-harcos} show that Bianchi newforms with Galois stable Nebentypus have associated 2-dimensional irreducible $p$-adic Galois representations. More precisely, let $f$ be a Bianchi newform\footnote{The term {\em newform} is used, both in elliptic and Bianchi setting, in its usual sense, that is, cuspidal, primitive and new eigenform.} of level 
$\mathcal{K}_0(\mathfrak{n})$\footnote{The group $\mathcal{K}_0(\mathfrak{n})$ is defined as 
$\left \{ \left ( \begin{smallmatrix} a & b \\ c&d \end{smallmatrix} \right )  \in \gl(\widehat{\O_K}) \mid c  \in \mathfrak{n}\widehat{\O_K} \right \}$. Notice that $\gl(\O_K) \cap \mathcal{K}_0(\mathfrak{n})$ is the standard congruence subgroup $\Gamma_0(\mathfrak{n})$.} with Nebentypus $\chi$ over an imaginary quadratic field $K$. It is essential for the techniques of Taylor et al. 
that we assume that $\chi=\chi^{c}$ where $c$ is the non-trivial automorphism of $K$. Let $p$ be a prime that is unramified in $K$ such that 
$(p)$ is coprime to $\mathfrak{N}$. Then there is an irreducible continuous representation 
$$\rho_{f,p} : G_K \longrightarrow \gl(\overline{\Q}_p)$$
of the absolute Galois group $G_K$ of $K$ such that the associated $L$-functions agree, that is,
$$L(f,s) =^* L(\rho_{f,p},s)$$
where $=^*$ means that the local Euler factors {\em away from} $(p) (D_K) \mathfrak{n} \mathfrak{n}^c$ agree. See \cite[Theorem A]{hst} and 
\cite[Section 3]{dgp} for details. Local properties of the above $\rho_{f,p}$ have been further investigated by Jorza in \cite{jorza}.

Due to the nature of the proof of the above statement, we do not know whether these $p$-adic representations arise, as expected, from the et\'ale 
cohomology of some complex algebraic variety. Constructing such varieties is one of the biggest challenges in the theory.   
Perhaps, the first step in attacking this problem should be to consider, as we shall do in the next section, the case of weight 2 Bianchi modular forms and 
elliptic curves over imaginary quadratic fields. 


\section{Connections with Elliptic Curves} \label{elliptic}

Work of Shimura \cite{shimurabook}, Wiles \cite{wiles}, Taylor-Wiles \cite{taylor-wiles} and Breuil-Conrad-Diamond-Taylor \cite{breuil} provide us with a celebrated bijection, via the associated $L$-functions, between isogeny classes of elliptic curves over $\Q$ of fixed conductor $N$ and elliptic 
newforms for $\Gamma_0(N) \leq \sl(\Z)$ of weight 2 with rational eigenvalues. In this section, we will discuss the conjectural analogous bijection, again defined 
via the associated $L$-functions, over imaginary quadratic fields.

Perhaps the first time such a connection over imaginary quadratic fields was mentioned is 1970. In his article \cite{serre-70}, which appeared very soon after Weil's influential paper \cite{weil}, Serre defines Hecke operators on the first homology of a Bianchi group $\Gamma$ 
and makes the following remark, see \cite[p.613]{serre-70}; 
\begin{quotation}
``Pour certains sous-groupes $\Gamma$, les r\'esultats r\'ecents de Weil laissent penser 
que les valeurs propres en question sont \'etroitement li\'ees aux propri\'et\'es arithm\'etiques
des courbes elliptiques d\'efinies sur $K$; il serait tr\'es int\'eressant d'en avoir des
exemples explicites."
\end{quotation}
 
The first investigations in this direction were done by Grunewald and Mennicke \cite{grunewald-mennicke} in 1978. They tabulated all elliptic curves of prime conductor with norm $\leq 1000$ and all Bianchi newforms for $\Gamma_0(\p)$ of weight 2 of prime level $\p$ with norm $\leq 1000$ over imaginary quadratic fields with $d=1,2,3$.  In 1981, Cremona \cite{cremona-84} investigated this phenomenon more extensively working over fields with $d=1,2,3,7,11$ and considering all elliptic curves of conductor $\leq 500,300,500,200,200$ over these fields respectively. Cremona's students extended his methods to several higher class number fields in \cite{whitley, bygott, lingham}. 

The data collected in all the above works strongly support the existence of such a bijection, however with certain natural modifications. 
If $E$ is an elliptic curve over $K$ with CM by an order in $\O_K$, then $E$ cannot be modular by a Bianchi newform, instead it corresponds to an Eisenstein series, see \cite{egm-80}. It was observed in \cite{egm-80, cremona-84} that for certain Bianchi newforms $F$ of weight 2 and rational Hecke eigenvalues, there was no elliptic curve $E$ over $K$ which is modular by $F$. As first observed by Deligne (see \cite{egm-80} p.267) and studied in detail by Cremona in \cite{cremona-92}, in such situations, it is expected that there is a {\em fake elliptic curve} that is modular by $f$, in other words, there is an abelian surface $A$ over $\Q$ such that $K \simeq End_{\Q}(A)$ and $End_{K}(A)$ is a division algebra. Such surfaces are modular by elliptic newforms $f$ with {\em extra twist}, see Cremona \cite[Theorem 5]{cremona-92}. Hence our original Bianchi newform $F$ is a twist by a quadratic character of $G_K$ of the base-change of $f$ to $K$. 

Keeping the above in mind, we end up with the following conjectures, see \cite{egm-80,cremona-92,dgp}. For general number fields, see Taylor 
\cite{taylor-95}.
\begin{svgraybox}
\begin{conjecture} \label{taniyama-shimura} Let $K$ be an imaginary quadratic field. 
\begin{enumerate}
\item Let $F$ be a Bianchi newform of level $\mathcal{K}_0(\mathfrak{n})$ and weight 2 over $K$ with rational eigenvalues. 
Assume that $F$ is not a twist by a quadratic character of $G_K$ of the base-change to $K$ of an elliptic newform. Then there is an elliptic curve
over $K$ of conductor $\a$ which is modular by $F$.   
\item Let $E$ be an elliptic curve over $K$ of conductor $\a$ which does not have CM by an order in $K$. Then $E$ is modular by some Bianchi newform $F$ of level 
$\mathcal{K}_0(\mathfrak{n})$ and weight 2 over $K$.
\end{enumerate}
\end{conjecture}
\end{svgraybox}

Given a Bianchi newform $F$ over $K$ as above and an elliptic curve over $K$ which, by virtue of numerical evidence, seems modular by $F$, it is in fact possible to {\em prove} that this is so via the so called {\em Faltings-Serre method}, see \cite{serre-84}. This is a method that allows one to decide, from a comparison of finitely many Euler factors of the associated $L$-functions, whether two given semisimple continuous representations of $G_K$ into $\gl(\overline{\Q}_p)$ are isomorphic or not. Hence, to be able to use this method for our purpose, the existence of Galois representations associated to Bianchi newforms is essential (this was discussed in the previous section). The first application 
of this method in our setting was carried out by Taylor at the end of \cite{taylor-94} where it was shown, using data provided by Cremona, that the elliptic curve 
over $\Q(\sqrt{-3})$ given by 
$$y^2 + x y = x^3 + \frac{3 + \sqrt{-3}}{2} \ x^2 + \frac{1+\sqrt{-3}}{2} \ x$$
is modular. This is a non-CM elliptic curve of conductor $\left ( \frac{17+\sqrt{-3}}{2} \right )$, with norm $73$, and is not isogenous to a curve that is defined over $\Q$. Recently, Dieulefait, Guerberoff and Pacetti developed \cite{dgp} an algorithm, based on the Faltings-Serre method, that 
outputs a finite set $S$ of rational primes such that a comparison of the Euler factors indexed by the prime ideals over $S$ allows one to decide whether a given elliptic curve over an imaginary quadratic field $K$ is modular or not by a given Bianchi newform of weight 2 over $K$. See \cite{dgp} for more examples of modular elliptic curves over imaginary quadratic fields. Detailed discussions of the Faltings-Serre method can be found in \cite{livne, chenevert, schutt}.

One of the most fundamental problems in the arithmetic theory of Bianchi groups is the construction of the elliptic curves that are 
conjectured to be associated to Bianchi newforms as above. 

\begin{svgraybox}
\begin{question} Given a Bianchi newform as in Conjecture \ref{taniyama-shimura}, how can we {\em geometrically} construct the conjecturally associated
elliptic curve ?
\end{question}
\end{svgraybox}
 
The principal obstruction to imitating Shimura's classical geometric method is the fact that the hyperbolic 3-folds $Y$ associated to Bianchi groups 
and their congruence subgroups lack complex structure and thus one does not have any (complex) algebraic variety structure on $Y$. In \cite{egm-80}, Elstrodt, Grunewald and Mennicke tried, however unsuccessfully, to obtain these elliptic curves from surfaces embedded in $Y$. The periods attached to a Bianchi newform of weight 2 with rational Hecke eigenvalues are all real and they form a one-dimensional lattice in $\R$ (see \cite{kurcanov}). This makes the well-known methods of constructing the elliptic curve from its period lattice unsuccessful.

\runinhead{Elliptic Curves Over Imaginary Quadratic Fields}
Elliptic curves over imaginary quadratic fields $K$ with restricted ramification properties have been studied by several authors, see 
\cite{comalada-nart, comalada, ishii-79, ishii-86, kida-01, pinch-84, pinch-86, pinch-87, setzer-78, setzer-81, stroeker}. Note that it is possible to list algorithmically all elliptic curves with prescribed ramification over a number field, see \cite{cremona-lingham, nook}. Especially interesting 
is the case of elliptic curves over $K$ with good reduction everywhere. It was shown by Stroeker \cite{stroeker} that such curves cannot have globally minimal 
models. Combined with an earlier result of Setzer \cite{setzer-78}, one concludes that imaginary quadratic fields whose class numbers are prime to 6 do not admit elliptic curves with good reduction everywhere. For a comprehensive list of imaginary quadratic fields which are known not to admit elliptic curves of 
good reduction everywhere, see the 2011 Ph.D. thesis of Thongjunthug \cite{nook}. 

It is a fact that goes back to Tate \cite{tate}, see also \cite{stroeker}, that there are infinitely many imaginary quadratic fields which admit elliptic curves with 
good reduction everywhere. In his 1984 Ph.D. thesis \cite{kraemer}, Kr\"amer proved that there are only eight such imaginary quadratic fields $K$ which 
have odd class number (see Table \ref{table: kraemer}). Kr\"amer also proved that the smallest, with respect to the absolute value of the discriminant, imaginary quadratic field which admits an elliptic curve with good reduction everywhere is $\Q(\sqrt{-259})$ (of class number $4$). By Setzer \cite{setzer-78}, we know that the next smallest is $\Q(\sqrt{-65})$ (of class number $8$). In \cite{scheutzow}, Scheutzow computed the dimensions of spaces of cuspidal Bianchi newforms of weight 2 and level one over $K$ with discriminant $\geq -260$. The data he collected, together with Kr\"amer's results, lend strong support to Conjecture \label{taniyama-shimura} above. 
 
\begin{table}
\centering 
\begin{tabular}{|c|c|} \hline
$d$ & $a \in \Q(\omega)$ \\ \hline
$643$ & $(1/2)(63+\omega)$ \\ \hline
$1879$ & $135+2\omega$ \\ \hline
$10691$ & $(1/2)(55+\omega)$ \\ \hline
$127285811$ & $(1/2)(-8441+\omega)$ \\ \hline
$218017403$ & $(1/2)(13427+\omega)$ \\ \hline
$349633643$ & $(1/2)(21867+\omega)$ \\ \hline
$443038171$ & $(1/2)(32355 + \omega)$ \\ \hline
$451767203$ & $(1/2)(-35297+\omega)$ \\ \hline
\multicolumn{2}{|c|}{$\omega = \sqrt{-d}$} \\
\multicolumn{2}{|c|}{$j=j(E)=a^2+1728$} \\ 
\multicolumn{2}{|c|}{$E: \ \ y^2=x^3-27j^3 \ x + 54aj^4$} \\ \hline 
\end{tabular}
\caption{All the odd class number imaginary quadratic fields which admit elliptic curves with everywhere good reduction. 
The unique, up to isomorphism, such elliptic curve is also described.} 
\label{table: kraemer}
\end{table}

It is a good place to bring to the attention of the reader a perhaps subtle issue. An elliptic curve $E$ over an imaginary quadratic field $K$ which has good reduction everywhere cannot be defined over $\Q$, see \cite{kraemer, kida-02}. However, the conjecturally associated Bianchi newform of weight 2 can be a lift of an elliptic newform of weight 2. Let us illustrate this through an explicit example. It is well-known that there is an abelian surface $A$ over $\Q$ of conductor $37^2$ with $End_\Q(A) \simeq \Q(i)$. This surface is associated to a pair of Galois conjugate elliptic newforms of weight 2 for $\Gamma_0(37)$ with Nebentypus $\chi_{37}$, the quadratic character given by the field $\Q(\sqrt{37})$. These two elliptic newforms lift to the unique level one weight 2 Bianchi newform $F$ over $K=\Q(\sqrt{-259})$. The base change of the surface $A$ to $K$ splits into two isomorphic elliptic curves with good reduction everywhere defined over $K$. These elliptic curves are modular by $F$.   

\runinhead{Higher Dimensional Abelian Varieties}
The above conjectural picture naturally extends to abelian varieties over $K$ with special types of endomorphism rings. An abelian variety $A$ of dimension $n$ over a number field is said to be {\em of GL(2)-type} if there is a number field $L$ of degree $n$ such that $End_K(A) \simeq L$. In this case, the action of the absolute Galois group $G_K$ of $K$ on the Tate modules of $A$ breaks into 2-dimensional representations. This leads to a notion of modularity, defined similarly by stating that the $L$-function of the abelian variety is given by the $L$-function of newform of weight $2$ and its Galois conjugates. See Ribet \cite{ribet} for details. It is a fact that every abelian variety $A$ over $\Q$ that is of $GL(2)$-type is modular. This follows from (see \cite[Theorem 4.4.]{ribet}) the celebrated Serre's Modularity Conjecture \cite{serre-87} which was recently proven by Khare and Wintenberger, see \cite{khare-wintenberger}. It is expected that abelian varities of $GL(2)$-type over any number field should be modular, see Taylor \cite{taylor-95}.

Compared to elliptic curves, higher dimensional abelian varieties of $GL(2)$-type are much harder to produce. There are several families of abelian surfaces (described via the Weierstrass equation of the defining hyperelliptic curve) with prescribed real multiplication in the literature, see Mestre \cite{mestre}, Brumer \cite{brumer} (contains Mestre's family), Bending \cite{bending} for multiplication by $\sqrt{2}$ and Wilson \cite{wilson} for multiplication by $\sqrt{5}$. The general problem with using the above families is that during the process of specialization to a number field $K$, one cannot control the ramification locus of the resulting abelian surface over $K$. Recently Elkies and Kumar \cite{elkies-kumar} developed an algorithmic method that produces a {\em universal} family of abelian surfaces (described via the Igusa-Clebsch invariants) with any prescribed real multiplication.  

A couple of years ago, A.Brumer and the author found several examples of abelian surfaces with real multiplication over imaginary quadratic fields which numerically seemed modular. Let us present one such example here. Let $\mathcal{C}$ be the hyperelliptic curve of genus $2$ given by 
$$y^2 = x^6+(8+2i)x^5+(25+10i)x^4+(40+14i)x^3+36x^2+(14-10i)x+(1-4i)$$
over $\Q(i).$ Then the Jacobian $A =Jac(\mathcal{C})$ of $\mathcal{C}$ is an abelian surface with multiplication by $\sqrt{5}$ whose conductor is $\mathfrak{N}^2$ where $\mathfrak{N}=(22+15i)$. There is a pair of Galois conjugate Hecke eigenclasses $f,f^{\sigma} \in H^1_{cusp}(\Gamma_0(\mathfrak{N}),\mathbb{C})$ with coefficient fields $\Q(\sqrt{5})$ for which our data suggest that
$$L(A,s)=L(f,s) \times L(f^{\sigma},s).$$


\section{The Torsion Phenomenon}  \label{torsion}
It is well-known that the $p$-torsion elements of an infinite group $\Gamma$ give rise to $p$-torsion classes in the (co)homology 
of $\Gamma$. For some groups $\Gamma$, one observes $p$-torsion in the (co)homology even when $\Gamma$ has no 
$p$-torsion. Such torsion was labeled ``{\em strange torsion}" by Weintraub in \cite{weintraub}. We will use our own terminology. 

Let $\Gamma$ be a finite index subgroup of a Bianchi group $\psl(\O)$. For $p \not= 2,3$, a $p$-torsion class in $H^i(\Gamma, V_{k,\ell}(\O))$ or 
$H_i(\Gamma, V_{k,\ell}(\O))$ is called {\bf large}\footnote{Note that if $p \leq max\{ k,\ell \}$, then $V_{k,\ell}(\O / (p))$ is a reducible 
$\psl(\O / (p))$-module and its coinvariants typically give rise to $p$-torsion in the homology.} if $p > max\{ k,\ell \}$. 
Existence of large torsion in the (co)homology raises interesting number theoretical questions as we shall discuss later.

Large torsion was first observed by Grunewald in \cite{egm-81} where he computed that for many small level congruence subgroups $\Gamma$ of Bianchi groups $\psl(\O_d)$ with $d=1,2,3$, the abelian group $H_1(\Gamma,\Z) \simeq \Gamma^{ab}$ has $p$-torsion with $p>3$. Note that Bianchi groups have only 2 and 3 torsion. Taylor in his 1988 Ph.D. thesis \cite{taylor-88}, for principal congruence subgroups $\Gamma$ with
large enough $p$-power level, showed the existence of $p$-torsion in $H_1(\Gamma,V_{k,\ell}(\O))$ for $k \not=\ell$. In 2000, Priplata \cite{priplata} found examples of torsion in $H_1(\psl(\O_d),V_{k,\ell}(\O_d))$ for $d=1,2,3$ and small $k \not= \ell$. Recently the author extended these 
computations in \cite{sengun-torsion} to $H_1(\Gamma,\Z)$ for very large level congruence subgroups $\Gamma$ and to parallel weights 
$H^2(\pgl(\O_d), V_{k,k}(\O_d))$ for $d=1,2,3,7,11$. Data showed that the $p$-torsion 
gets to astronomical sizes. We give a modest sample in Tables \ref{torsion1} and \ref{torsion2}, see \cite{sengun-torsion, rahm-abelianization} for more.

\begin{table}
\centering 
\begin{tabular}{|c|c|} \hline
$k$ & $p$ \\ \hline
24&  2, 3, 5, 7, 11, 13, 17, 23, {\bf 1033}  \\ \hline   
25&  2, 3, 5, 7, 11, 17  \\ \hline     
26&  2, 3, 5, 7, 11, 13, 19, 23, {\bf 157}, {\bf 683}  \\ \hline 
27&  2, 3, 5, 7   \\ \hline 
28&  2, 3, 5, 7, 11, 13, 17, {\bf 664197637}      \\ \hline      
29&  2, 3, 5, 7, 11, 13, {\bf 89}          \\ \hline           
30&  2, 3, 5, 7, 11, 13, 19, 23, 29, {\bf 211}, {\bf 36312691}   \\ \hline  
\end{tabular}
\caption{All the $p$-torsion in $H^2(\pgl(\Z[i]),E_{k,k}(\Z[i]))$. The ones in bold are large torsion.}
\label{torsion1}
\end{table}

\begin{table}
\centering 
\begin{tabular}{|c|c|} \hline
$\textbf{N}\mathfrak{p}$ &  $p$ \\ \hline
4909 &  7, 13, 409, 10691, 22871, 29423, 56980673, 71143433 ] \\ \hline
4933 & 37, 101, 137, 577, 947, 21169, 194981 ] \\ \hline
4937 &  7, 37, 617, 10859, 108893, 4408530403, 157824962047 ] \\ \hline
4957 &  7, 13, 31, 59, 14347, 3051863, 9405667, 23132267 ] \\ \hline
4969 &  23, 71, 373, 191299, 39861006443, 8672729371087 ] \\ \hline
4973 &  11, 13, 47, 71, 113, 127, 331, 6317, 7949, 39023, 628801, 2995319 ] \\ \hline
4993 &  5, 7, 11, 13, 101, 173, 798569, 5995036891, 18513420749 ] \\ \hline
\end{tabular}
\caption{{\em Some} of the large $p$-torsion in $H_1(\Gamma_0(\p),\Z)$ for $\Gamma_0(\p) \subset \psl(\Z[i])$}
\label{torsion2}
\end{table}

Recent works of Bergeron-Venkatesh \cite{bergeron-venkatesh} and Marshall-M\"uller \cite{marshall-muller} brought a very interesting perspective to the large torsion phenomenon. In \cite{bergeron-venkatesh}, Bergeron and Venkatesh studied the growth of torsion in the homology of cocompact arithmetic groups. Their strongest result is in the case where the ambient Lie group is $\psl(\C)$. 
Let $\Gamma$ be a cocompact arithmetic lattice in $\psl(\C)$. Let  $R$ be the 
ring obtained from the ring of integers of the invariant trace field of $\Gamma$ by inverting the minimal set of
of prime ideals which ensures that $\Gamma \subset \psl(R)$. 
For a pair of non-negative integers $k$ and $\ell$ define
$$c_{k,\ell}:=\dfrac{1}{6\pi}\cdot \dfrac{1}{2^3} \cdot
\Big (  (k+\ell+2)^3 - |k-\ell |^3 + 3(k+\ell+2)(k-\ell)(k+\ell+2-|k-\ell|) \Big ).$$
Note that $c_{0,0}=\frac{1}{6\pi}.$

Let $\{ \Gamma_n \}_{n \geq 1}$ be a decreasing tower of congruence subgroups of $\Gamma$
such that $\bigcap_n \Gamma_n = \{ 1 \}$. Then for $k \not= \ell$, Bergeron and Venkatesh prove 
that
$$ \lim_{n \rightarrow \infty} \dfrac{\log | H_1(\Gamma_n, V_{k,\ell}(R))_{tor} |}{\textrm{vol}(\Gamma_n \backslash \H)}
=c_{k,\ell}.$$

As they are not cocompact, Bianchi groups and their congruence subgroups are outside the scope of the above result. Moreover, 
the above result only considers acyclic modules, that is, $V_{k,\ell}$ with $k \not= \ell$. As we have discussed above, from a number 
theoretical perspective, the modules $V_{k,k}$ are more interesting. In \cite{sengun-torsion}, the author collected numerical data 
which strongly suggest that the above result of Bergeron and Venkatesh extend to non-cocompact arithmetic lattices 
in $\psl(\C)$, to non-acyclic modules and to collections of congruence subgroups which satisfy weaker conditions 
(see \cite{7samurai, 7samurai-paper} for such a condition) as well. In particular, we have the following.

\begin{svgraybox}
 \begin{conjecture} Consider congruence
subgroups of the form $\Gamma_0(\p)$ inside a fixed Bianchi group. Then
$$ \lim_{N\p \rightarrow \infty} \dfrac{\log | H_1(\Gamma_0(\p), \Z)_{tor} |}{\textrm{vol}(\Gamma_0(\p) \backslash \H)}
=\dfrac{1}{6\pi}$$
where the limit is taken over prime ideals $\p$ of residue degree one.
\end{conjecture}
\end{svgraybox}

For a more general conjecture, see Section 10 of the recent paper \cite{7samurai-paper} by Abert et al.

\section{Connections with Galois Representations II} \label{galois2}
In this section we will discuss connections between eigensystems occuring in the mod $p$ homology of 
congruence subgroups $\Gamma$ of Bianchi groups $\pgl(\O_K)$ and continuous representations 
of the absolute Galois group $G_K$ of $K$ into $\textrm{GL}_2(\overline{\F}_p)$. 
This connection is in the sense of Serre \cite{serre-87} and Ash \cite{ash-89}.

Fix a prime ideal of $\O_K$ over the rational prime $p$ with residue field $\kappa$. Let us analyze the following diagram
$$\xymatrix{ H^2(\Gamma,V_{k,\ell}(\mathcal{O}_K)) \ar[r]^{\otimes \C}  \ar@{->>}[d]^{\otimes \kappa } & 
H^2(\Gamma,V_{k,\ell}(\C)) \\  H^2(\Gamma,V_{k,\ell}(\overline{\F}_p)) }$$
from the perspective of the Hecke action. It is well-known that the image of the integral cohomology 
inside the complex cohomology forms a basis and thus all the arithmetic information that the complex cohomology carries is captured 
by the torsion-free part of the integral cohomology. On the other hand, for large $p$, the arithmetic information captured by 
the mod $p$ cohomology lifts to the integral cohomology (see \cite{ash-stevens, sengun-torsion}). The critical issue is that 
due to the existence of possible large $p$-torsion in the integral cohomology, there may be eigenvalue systems in the mod 
$p$ cohomology which are not reductions of eigenvalue systems captured in the complex cohomology and therefore 
they are not, a priori, related to automorphic forms. Yet it is believed that such mod $p$ eigenvalue systems have associated 
2-dimensional mod $p$ Galois representations as well, see Ash \cite{ash-89}. Conversely, a natural generalization of Serre's 
(weak) modularity conjecture predicts that any continuous semi-simple representation 
$G_K \rightarrow \gl(\overline{\F}_p)$ should be associated to a mod $p$ eigenvalue system. Note that 
the ``odd/even" distinction that exists for representations of $G_\Q$ does not exist for representations 
of $G_K$. 

The first numerical examples related to the above were found by Grunewald in \cite{egm-81} where the data collected 
suggested that several representations $G_{\Q(i)} \rightarrow \textrm{GL}_2(\F_2)$ were associated to eigenvalue 
systems in the first homology of certain congruence subgroups of $\psl(\Z[i])$ with coefficients in $\F_2$. 
Figueiredo studied \cite{figueiredo-98, figueiredo-99} Serre's modularity conjecture over the 
imaginary quadratic fields with $d=1,2,3$. For these fields, he found three representations $G_K \rightarrow\textrm{GL}_2(\F_3)$ 
which seemed to be associated to eigenvalue systems in the first homology of suitable congruence subgroups with coefficients 
in $\overline{\F}_3$. Figueiredo remarks that the mod $3$ eigenvalue systems that came up in his computations were not reductions of 
eigenvalue systems that were captured in the complex homology. See \cite{sengun-numerical} for further examples of this kind. 

Recently in \cite{torrey}, Torrey investigated a strong version 
of Serre's modularity conjecture over imaginary quadratic fields using the weight recipe of Buzzard, Diamond and Jarvis 
\cite{bdj} which was formulated for the Hilbert modular groups. She studied several numerical examples, including those of 
Figueiredo, from the perspective of weights. 

\subsection{Even Mod $p$ Galois Representations} \label{even}
Even though Serre's modularity conjecture for $\Q$ necessarily leaves out the {\em even} 
representations of $G_\Q$, we can talk about the modularity of such representations using 
Serre's modularity conjecture over imaginary quadratic fields as follows.

Start with an irreducible even continuous representation $\rho: G_\Q \rightarrow \textrm{GL}_2(\overline{\F}_p)$. 
Let $\rho_K $ be its restriction to $G_K$ for some imaginary quadratic field $K$. By Serre's modularity conjecture over $K$, $\rho_K$ 
should be associated to an eigenvalue system $\Psi$ occuring in $H^2(\Gamma, V(\overline{\F}_p))$ for some suitable congruence subgroup of 
$\Gamma$ of $\textrm{GL}_2(\O_K)$ and some irreducible $\F_p[\gl(\O_K / (p))]$-module $V(\overline{\F}_p)$. In fact, 
the module $V(\overline{\F}_p)$ can be specified using Serre's original weight recipe, or more generally, using the Buzzard-Diamond-Jarvis 
weight recipe, as done by Torrey. 

Let us illustrate the above through a concrete example. Consider the polynomial  
$$x^8 - x^7 - 44x^6 + 43x^5 + 442x^4 - 32x^3 - 1311x^2 - 1156x - 241.$$ 
Its Galois group is $\widehat{A_4}$ which has an embedding into $\textrm{GL}_2(\Z[\zeta_3])$.
The arising Artin representation
$$\rho : G_{\Q} \rightarrow \textrm{GL}_2(\C)$$ 
is even and irreducible. Moreover it has Artin conductor $163^2$. 

We restrict $\rho$ to $G_{\Q(i)}$ and then reduce the (finite) image mod $163$ to get 
$$\bar{\rho}_K : G_{\Q(i)} \rightarrow \textrm{GL}_2(\F_{163})$$
with trivial Serre conductor.

Using the computer programs written by the author to compute the Hecke action on the cohomology of 
the Euclidean Bianchi groups, we found  a seemingly matching eigenvalue system $\Psi$ occuring in the mod $163$ 
cohomology space
$$H^2(\textrm{GL}_2(\Z[i]),V(\F_{163}))$$
for weights
$$V_{53,53}^{27,27}(\F_{163}) \ \textrm{and} \ V_{107,107}^{135,135}(\F_{163}).$$ 
Here 
$$V_{k,\ell}^{a,b}(\F_{163}):= \left ( V_k(\kappa)\otimes_{\kappa} \textrm{det}^a \right ) 
\otimes_{\kappa} \left ( V_{\ell}(\kappa)^{Frob}\otimes_{\kappa} \textrm{det}^b \right )$$
with $\kappa \simeq \F_{163^2}$ being the residue field of  the prime ideal $(163)$ of $\Z[i]$ generated by 163. The group $\textrm{GL}_2(\Z[i])$ acts on $V_{k,\ell}^{a,b}$ via reduction modulo  $(163)$.

It turns out, see  \cite{sengun-even}, that $\Psi$ is {\em not} the mod $163$ reduction of an eigenvalue system occuring in the complex cohomology space  
$$ H^2(\textrm{GL}_2(\Z[i]),V(\C))$$ 
for the weights $V_{53,53}^{27,27}$ and $V_{107,107}^{135,135}$.

We conclude the existence of a 163-torsion class which is repsonsible for the eigenvalue system $\Psi$. 
In fact, we verified by direct computation as well that $H^2(\textrm{GL}_2(\Z[i]),V(\Z[i]))$ has 163-torsion 
for the weights $V_{107,107}^{135,135}$ and $V_{53,53}^{27,27}$. Moreover, we checked that all of the above 
remain true when we replace $\Q(i)$ by the other Euclidean imaginary quadratic fields. 

We speculate that for any irreducible, even, continuous representation $\rho: G_\Q \rightarrow \gl(\overline{\F}_p)$ and 
any imaginary quadratic field $K$, the eigenvalue system $\Psi$ that we will capture, 
following the above steps, will not be the reduction of an eigenvalue system that occurs in the relevant complex cohomology space. 
To motivate our speculation, let us assume that $\rho$, as in our example above, has trivial Serre conductor and that $\Psi$ 
is the mod $p$ reduction of an eigenvalue system captured by a level one cuspidal Bianchi modular form $F$. 
Then our discussion before Question \ref{bc} suggests that $F$ is most likely the base change to $K$ of 
an elliptic newform $f$. It follows that the {\em odd} representation $\rho_f : G_\Q \rightarrow \gl(\overline{\F}_p)$ 
associated to $f$ and our even representation $\rho$ become isomorphic when restricted to $G_K$. Basic Clifford theory shows that this is not 
possible.

\begin{acknowledgement}
It is a pleasure to thank Lassina Demb\'el\'e for his help with the part about connections with elliptic curves. We also thank T.Berger, J.-P. Serre and 
A.Rahm for their helpful feedback. Figures 1 and 2 were reproduced from Maite Aranes' thesis with her kind permission. Finally, we thank the 
referee for the constructive comments which gave the paper its more coherent final form.
\end{acknowledgement}

\bibliographystyle{siam}
\bibliography{SurveyReferences}
\end{document}